\documentclass[review]{elsarticle}
\usepackage{amsmath,amssymb}
\usepackage{multirow}
\usepackage{enumerate}
\usepackage{array}
\usepackage{algorithm}
\usepackage{algorithmic}
\usepackage{xcolor}

\newtheorem{thm}{Theorem}

\newtheorem{lmm}[thm]{Lemma}
\newtheorem{prop}[thm]{Proposition}
\newproof{pf}{Proof}

\journal{European Journal of Combinatorics}

\begin{document}
\begin{frontmatter}
\title{Enumeration and classification of self-orthogonal partial Latin rectangles by using the polynomial method.}
\author{Ra\'ul M. Falc\'on}
\ead{rafalgan@us.es}
\cortext[cor1]{Corresponding author. Phone: + 34 954 550158; Fax: +34 954 556683}

\address{School of Building Engineering, University of Seville,\\ Avenida de Reina Mercedes 4 A, 41012, Seville, Spain.}

\begin{abstract}
The current paper deals with the enumeration and classification of the set $\mathcal{SOR}_{r,n}$ of self-orthogonal $r\times r$ partial Latin rectangles based on $n$ symbols. These combinatorial objects are identified with the independent sets of a Hamming graph and with the zeros of a radical zero-dimensional ideal of polynomials, whose reduced Gr\"obner basis and Hilbert series can be computed to determine explicitly the set $\mathcal{SOR}_{r,n}$. In particular, the cardinality of this set is shown for $r\leq 4$ and $n\leq 9$ and several formulas on the cardinality of $\mathcal{SOR}_{r,n}$ are exposed, for $r\leq 3$. The distribution of $r\times s$ partial Latin rectangles based on $n$ symbols according to their size is also obtained, for all $r,s,n\leq 4$.
\end{abstract}

\end{frontmatter}

\section{Introduction.}

An {\em $r\times s$ partial Latin rectangle} based on $[n]=\{1,\ldots,n\}$ is an $r \times s$ array in which each cell is either empty or contains a symbol of $[n]$, such that each symbol occurs at most once in each row and in each column. Its number of filled cells is its {\em size}. Let $\mathcal{R}_{r,s,n}$ and $\mathcal{R}_{r,s,n:m}$ respectively denote the set of $r\times s$ partial Latin rectangles based on $[n]$ and its subset of partial Latin rectangles of size $m$. Given $P=(p_{ij})\in \mathcal{R}_{r,s,n}$, its {\em orthogonal array representation} is the set $O(P)=\{(i,j,p_{ij})\colon\, i\in [r], j\in [s], p_{ij}\in [n]\}$. Permutations of rows, columns and symbols of $P$ give rise to new $r\times s$ partial Latin rectangles based on $[n]$, which are said to be {\em isotopic} to $P$. If $S_m$ denotes the symmetric group on $m$ elements, then $\Theta=(\alpha,\beta,\gamma)\in S_r\times S_s\times S_n$ is an {\em isotopism} of $\mathcal{R}_{r,s,n}$ and it is defined the isotopic partial Latin rectangle $P^{\Theta}$ such that $O(P^{\Theta})=\{(\alpha(i),\beta(j),\gamma(p_{i,j}))\colon\, i\in [r], j\in [s], p_{ij}\in [n]\}$. Given a permutation $\pi\in S_3$, it is defined the {\em parastrophic} partial Latin rectangle $P^{\pi}$ such that $O(P^{\pi})=\{(p_{\pi(1)},p_{\pi(2)},p_{\pi(3)})\colon\, (p_1,p_2,p_3)\in O(P)\}$. If $P^{\pi}\in \mathcal{R}_{r,s,n}$, then $\pi$ is called a {\em parastrophism} of $\mathcal{R}_{r,s,n}$. The composition of an isotopism and a parastrophism is a {\em paratopism}. Two partial Latin rectangles are in the same {\em main class} if one of them is isotopic to a parastrophic partial Latin rectangle of the other.

Two partial Latin rectangles $P=(p_{ij}), Q=(q_{ij})\in \mathcal{R}_{r,s,n}$ are {\em orthogonal} if, given $i,i'\in [r]$ and $j,j'\in [s]$ such that $p_{ij}=p_{i'j'}\in [n]$, then $q_{ij}$ and $q_{i'j'}$ are not the same symbol of $[n]$. If $r=s$, then the partial Latin rectangle $P$ is {\em self-orthogonal} if it is orthogonal to its transpose $P^t$ (see Figure \ref{fig_PLR}).

\begin{figure}[h]
{\scriptsize
$$\begin{array}{|c|c|c|c|c|}\hline
           1 & 3 & \ \ & \ \ \\ \hline
           2 & \ \ & 3 & 1 \\ \hline
           \ \ & 1 & 2 & \ \ \\ \hline
           \ \ & 2 & \ \ & 3 \\ \hline
         \end{array}$$}
\caption{Example of a self-orthogonal $4\times 4$ partial Latin rectangle based on $[3]$.}
\label{fig_PLR}
\end{figure}

Let $\mathcal{SOR}_{r,n}$ be the set of self-orthogonal $r\times r$ partial Latin rectangles based on $[n]$. Only those isotopisms of the form $(\alpha,\alpha,\gamma)\in S_r\times S_r\times S_n$ and those paratopisms based on $\overline{S}_3=\{(1)(2)(3),(12)(3)\}$ preserve always the set $\mathcal{SOR}_{r,n}$. Hence, the sets $S_r\times S_n$ and $S_r\times S_n\rtimes \overline{S}_3$ determine, respectively, the {\em isotopism} and {\em paratopism groups} of $\mathcal{SOR}_{r,n}$. The enumeration of isotopism and main classes of $\mathcal{SOR}_{r,n}$ has been studied for $r=n\leq 10$ \cite{Burger2010, Burger2010a, Graham2006}. However, there does not exist a similar study for self-orthogonal partial Latin rectangles of any order. In the current paper, we deal with this problem by adapting the {\em Combinatorial Nullstellensatz} of Alon \cite{Alon1995}, whose effectiveness in the study of Latin squares has been exposed in \cite{Falcon2013,Falcon2007}.

The paper is organized as follows. In Section 2, we indicate some preliminaries concepts and results on commutative algebra. In Section 3, the set $\mathcal{R}_{r,s,n}$ is identified with that of independent sets of a Hamming graph and with the set of zeros of a zero-dimensional radical ideal, whose reduced Gr\"obner bases and Hilbert series determine, respectively, the elements and cardinality of $\mathcal{R}_{r,s,n:m}$, for all natural $m$. This cardinality is explicitly shown for $r\leq s\leq n\leq 4$ and $m\leq rsn$. In Section 4, we consider new polynomials to be added to the above ideal in order to determine the set $\mathcal{SOR}_{r,n}$. Besides, two strategies are indicated that allow us to reduce the cost of computation of the Gr\"obner basis and Hilbert series of the new ideal. They are used to determine the cardinality of $\mathcal{SOR}_{r,n}$ for $r\leq 4$ and $n\leq 9$. Some general formulas about the cardinality of $\mathcal{SOR}_{r,n}$ are finally exposed, for $r\leq 3$.

\section{Preliminaries.}

We start with some basic concepts of commutative algebra (see \cite{Cox1998, Cox2007, Kreuzer2005} for more details). Let $R=k[{\bf x}]=k[x_1,\ldots,x_n]$ be a polynomial ring in $n$ variables over a field $k$ with the {\em standard grading} induced by the degree of polynomials, that is, $R=\bigoplus_{0\leq d} R_d$, where each $R_d$ is the set of homogeneous polynomials in $R$ of degree $d$. The largest monomial of a polynomial of $R$ with respect to a given term order $<$ is its {\em initial monomial}. Given an ideal $I$ of $R$, the ideal generated by the initial monomials with respect to $<$ of all the non-zero elements of $I$ is its {\em initial ideal} $I_<$. Any monomial of $R$ which is not contained in $I_<$ is called a {\em standard monomial} of $I$ with respect to $<$. The set of standard monomials of $I$ with respect to any given term order can be used to study the dimension of the quotient ring $R/I$. This ring inherits the natural grading of $R$ and can be written as the direct sum $\bigoplus_{0\leq d} R_d/I_d$, where $I_d=R_d\cap I$. In particular, the set of standard monomials of $I$ of degree $d$ with respect to a given term order constitutes a linear $k$-basis of $R_d/I_d$ and hence, its cardinality coincides with $\mathrm{dim}_k(R_d/I_d)$, regardless of the term order which has been chosen. The {\em Hilbert function} $\mathrm{HF}_{R/I}$ of $R/I$ maps each non-negative integer $d$ onto $\mathrm{dim}_k(R_d/I_d)$. Its {\em Hilbert series} is the generating function $\mathrm{HS}_{R/I}(t)=\sum_{0\leq d}\mathrm{HF}_{R/I}(d)\cdot t^d$, which can also be written as:
\begin{equation}
\mathrm{HS}_{R/I}(t)=\frac {P(t)}{(1-t)^n}=\frac{Q(t)}{(1-t)^{\mathrm{dim}_k(I)}},
\end{equation}
where $P(t)$ and $Q(t)$ are polynomials with integer coefficients. The former is called the {\em Hilbert numerator} of $R/I$. If $I$ is zero-dimensional, then the Hilbert series of $R/I$ coincides with the polynomial $Q(t)$. The number of standard monomials of $I$ is then finite and coincides with the dimension of $R/I$ over $k$. Moreover, it is always greater than or equal to the number of points of the {\em affine variety} $V(I)$ of $I$, that is, the set of points in $k^n$ which are zeros of all the polynomials of $I$. The ideal $I$ is zero-dimensional if and only if $V(I)$ is finite. Further, it is verified that $|V(I)|=\mathrm{dim}_k(R/I)$ whenever $I$ is {\em radical}, that is, if any polynomial $p$ belongs to $I$ whenever there exists a natural $n\in \mathbb{N}$ such that $p^n\in I$.

The problem of computing a Hilbert series is NP-complete \cite{Bayer1992}. The standard algorithm which is commonly used to compute the Hilbert series of a quotient ring $R/I$ by determining the corresponding Hilbert numerator was first proposed by Mora and M\"oller in \cite{Mora1983} and is based on the additivity of the Hilbert function in short exact sequences. Even if it only works when $I$ is homogeneous, it is not an inconvenient, because the Hilbert series of $R/I$ coincides with that of $R/I_<$ for any term order $<$ and the initial ideal $I_<$ is homogeneous because all its elements are monomials. It is therefore necessary to determine explicitly the initial ideal of $I$ for a given term order $<$. In this regard, a {\em Gr\"obner basis} of $I$ with respect to $<$ is any generating set $G$ of $I$ such that the initial monomials of its elements generate the initial ideal $I_<$. It is said to be {\em reduced} if all its polynomials are monic and no monomial of a polynomial of $G$ can be generated by the initial monomials of the other polynomials of the basis. There exists only one reduced Gr\"obner basis of an ideal and the algorithm which is most commonly used to obtain it is that given by Buchberger \cite{Buchberger2006}.

\section{Boolean ideals and Hilbert series related to $\mathcal{R}_{r,s,n}$.}

In the current section, we show how the set $\mathcal{R}_{r,s,n}$ of $r\times s$ partial Latin rectangles based on $[n]$ can be identified with that of zeros of a boolean ideal which is zero-dimensional and radical. The use of non-linear polynomials to solve combinatorial problems was established by Alon \cite{Alon1995}. Afterwards, Bernasconi et al. \cite{Bernasconi1997} exposed the possibility of solving counting problems in Combinatorics by using Gr\"obner bases of {\em boolean} ideals, that is to say, ideals on $k[x_1,\ldots,x_n]$ containing the polynomials $x_i\cdot (x_i-1)$, for all $i\in [n]$. Such a boolean structure facilitates the computation of the corresponding reduced Gr\"obner basis with respect to any given term order.

\renewcommand{\baselinestretch}{1}
\begin{figure}[h]
\begin{center}
\begin{tabular}{cc}
\begin{tabular}{c}
$P\equiv\ $\begin{tabular}{|c|c|c|}\hline
1 &   & 4\\ \hline
  & 3 & 2\\ \hline
\end{tabular}\\ \\
$O(P)=\left\{\begin{array}{c}(1,1,1),(1,3,4),\\(2,2,3),(2,3,2)\end{array}\right\}$
\end{tabular}
&
\begin{tabular}{c}
\includegraphics[width=0.5\textwidth]{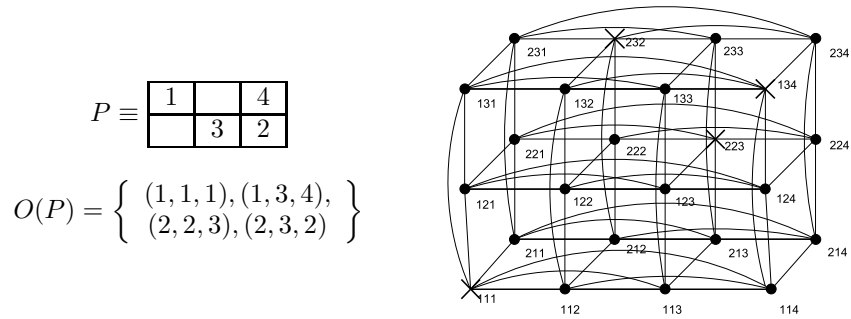}
\end{tabular}
\end{tabular}
\caption{$2\times 3$ partial Latin rectangle based on $[4]$ related to an independent set of $\mathcal{H}_{2,3,4}$.}
\label{fig0}
\end{center}
\end{figure}
\renewcommand{\baselinestretch}{2}

Let $\mathcal{H}_{r,s,n}$ be the {\em Hamming graph} \cite{Hammack2011} defined as the cartesian product $K_r\Box K_s\Box K_n$ of the three complete graphs of $r$, $s$ and $n$ vertices. It is a $(r+s+n-3)$-regular graph of order $rsn$, whose vertices can be labeled with the triples of $[r]\times [s]\times [n]$ so that two vertices are adjacent if and only if their corresponding labels differ exactly in one component. The labels of any {\em independent set} of $\mathcal{H}_{r,s,n}$, formed by $m$ pairwise nonadjacent vertices, constitute the orthogonal array representation of an $r\times s$ partial Latin rectangle based on $[n]$ of size $m$ (see Figure \ref{fig0}). The reciprocal is also verified and hence, the number of independent sets of $m$ vertices of $\mathcal{H}_{r,s,n}$ coincides with the cardinality of $\mathcal{R}_{r,s,n:m}$. When $r=s=n$, the Hamming graph $\mathcal{H}_{n,n,n}$ is usually denoted as $H(3,n)$ and its number of independent sets is only known for $n\leq 4$ (see the integer sequences A027681 and A027682 in \cite{Sloane2014}). The following result holds.

\begin{prop} \label{prop0} It is verified that:
\begin{itemize}
\item[a)] $|\mathcal{R}_{r,s,n:0}|=1$.
\item[b)] $|\mathcal{R}_{r,s,n:1}|=rsn$.
\item[c)] $|\mathcal{R}_{r,s,n:2}|=\frac 12 rsn\cdot (rsn-r-s-n+2)$.
\end{itemize}
\end{prop}

\begin{pf} The first two assertions are immediate. To prove the third one, observe that, since the number of independent sets of size two of any graph of $v$ vertices and $e$ edges is $\binom{v}{2} - e$, we have that $|\mathcal{R}_{r,s,n:2}|=\binom{rsn}{2}-\frac 12 \cdot rsn\cdot (r+s+n-3)=\frac 12 rsn\cdot (rsn-r-s-n+2)$. \hfill $\Box$
\end{pf}

\hspace{0.5cm}

Let us consider now a boolean variable $x_{ijk}$ related to each vertex of $\mathcal{H}_{r,s,n}$, for all $(i,j,k)\in [r]\times [s]\times [n]$, such that, given an independent set $S$ of $\mathcal{H}_{r,s,n}$, it takes the value $1$ if the vertex labeled as $(i,j,k)$ belongs to $S$, or $0$, otherwise. Keeping in mind the adjacency in our graph, the next result is satisfied.

\begin{thm} \label{thm0} The set $\mathcal{R}_{r,s,n}$ can be identified with the set of zeros of the ideal $I_{r,s,n}=\langle\,x_{ijk}(x_{ijk} -1),\, x_{ijk}x_{i'jk},\, x_{ijk}x_{ij'k},\, x_{ijk}x_{ijk'}\colon\, i\in [r],\, j\in [s],\, k\in [n], i'\in \{i+1,\ldots,r\}, j'\in \{j+1,\ldots,s\}, k'\in \{k+1,\ldots,n\}\,\rangle\, \subseteq \mathbb{F}_2[{\bf x}]$, where ${\bf x}=\{x_{111},\ldots,x_{rsn}\}$. Besides, $|\mathcal{R}_{r,s,n}|= \mathrm{dim}_{\mathbb{F}_2}(\mathbb{F}_2[{\bf x}]/I_{r,s,n})$ and  $|\mathcal{R}_{r,s,n:m}|=\mathrm{HF}_{\mathbb{F}_2[{\bf x}]/I_{r,s,n}}(m)$, for all $m\geq 0$.
\end{thm}

{\bf Proof.}  Any partial Latin rectangle $P=(p_{ij})\in\mathcal{R}_{r,s,n}$ can be uniquely identified with a zero $(x_{111},\ldots,x_{rsn})$, where $x_{ijk}=1$ if $p_{ij}=k$ and $0$, otherwise. The finiteness of $\mathcal{R}_{r,s,n}$ implies $I_{r,s,n}$ to be zero-dimensional. Further, for all $(i,j,k)\in [r]\times [s]\times [n]$, the intersection between $I_{r,s,n}$ and the polynomial ring $\mathbb{F}_2[x_{ijk}]$ of polynomials in the variable $x_{ijk}$ is the ideal $\langle\,x_{ijk}\left(x_{ijk}-1\right)\,\rangle\subseteq I_{r,s,n}$. Hence, Proposition 2.7 of \cite{Cox1998} assures $I_{r,s,n}$ to be radical and thus, Theorem 2.10 of \cite{Cox1998} implies that $|\mathcal{R}_{r,s,n}|=|V(I_{r,s,n})|= \mathrm{dim}_{\mathbb{F}_2}(\mathbb{F}_2[{\bf x}]/I_{r,s,n})$.

The reduced Gr\"obner basis of $I_{r,s,n}$ with respect to the lexicographic order $<_{\mathrm{lex}}$ coincides with its set of generators. Moreover, if we change in such a set each polynomial $x_{ijk}(x_{ijk}-1)$ by its initial monomial $x_{ijk}^2$ with respect to $<_{\mathrm{lex}}$, then the new set generates the initial ideal ${I_{r,s,n}}_{<_{\mathrm{lex}}}$. Let ${\bf x^a}=x_{111}^{a_{111}}\ldots x_{rsn}^{a_{rsn}}$ be a standard monomial of $I_{r,s,n}$ with respect to $<_{\mathrm{lex}}$. Since such a monomial does not belong to ${I_{r,s,n}}_{<_{\mathrm{lex}}}$, it is $a_{ijk}\in\{0,1\}$, for all $(i,j,k)\in [r]\times[s]\times[n]$. The standard monomial ${\bf x^a}$ can then be identified with the partial Latin rectangle $P=(p_{ij})\in\mathcal{R}_{r,s,n}$, such that $p_{ij}=k$ if $a_{ijk}=1$ and $\emptyset$, otherwise. The degree of ${\bf x^a}$ coincides with the size of $P$ and hence, given $m\geq 0$, the cardinality of $\mathcal{R}_{r,s,n:m}$ is equal to the number of standard monomials of degree $d$ of $I_{r,s,n}$. \hfill $\Box$

\vspace{0.5cm}

The initial ideal ${I_{r,s,n}}_{<_{\mathrm{lex}}}$ which appears in the proof of Theorem \ref{thm0} is the {\em modified edge ideal} \cite{Dickenstein2012} of the graph of $rsn$ vertices labeled by the variables $x_{111},\ldots,x_{rsn}$ and such that there exists an edge connecting the vertices labeled by $x_{ijk}$ and $x_{i'j'k'}$ if and only if the monomial $x_{ijk}x_{i'j'k'}$ belongs to ${I_{r,s,n}}_{<_{\mathrm{lex}}}$. Its standard monomials can then be identified with the independent sets of such a graph. It is the fundamental of a specialized algorithm exposed by Dickenstein and Tobis \cite{Dickenstein2012} to compute the Hilbert series related to this kind of ideals. We have implemented this algorithm in a procedure called {\em PLR} in the open computer algebra system for polynomial computations {\sc Singular} \cite{Decker2014}. It has been included in the library {\em pls.lib}, which is available online on {\tt http://personal.us.es/raufalgan/LS/pls.lib} and has been used to obtain the distribution in Table \ref{table_hs1} of the partial Latin rectangles of $\mathcal{R}_{r,s,n}$ with $r\leq s\leq n\leq 4$ according to their sizes. The computation of the corresponding Hilbert series is immediate ($0$ seconds) in an {\em Intel Core i7-2600, 3.4 GHz} with {\em Ubuntu}, except for $\mathcal{R}_{4,4,4}$, for which the CPU running time is $50$ seconds. Remark that all the values of Table \ref{table_hs1} have also been checked by exhaustive search.

\renewcommand{\baselinestretch}{1}
\renewcommand{\tabcolsep}{0.1pt}
\begin{table}[htb]
\begin{center}{\tiny
\hspace{-0.7cm} \begin{tabular}{rrrrrrrrrrrrrrrrrrrrrrr} \hline
\ & \multicolumn{5}{l}{$|\mathcal{R}_{r,s,n:m}|$}\\ \cline{2-21}
\ & \multicolumn{5}{l}{$r.s.n$}\\ \cline{2-21}
$m$ & 1.1.1 & 1.1.2 & 1.1.3 & 1.1.4 & 1.2.2 & 1.2.3 & 1.2.4 & 1.3.3 & 1.3.4 & 1.4.4 & 2.2.2 & 2.2.3 & 2.2.4 & 2.3.3 & 2.3.4 & 2.4.4 & 3.3.3 & 3.3.4 & 3.4.4 & 4.4.4\\ \hline
0 & 1 & 1 & 1 & 1 & 1 & 1 & 1 & 1 & 1 & 1 & 1 & 1 & 1 & 1 & 1 & 1 & 1 & 1 & 1 & 1\\
1 & 1 & 2 & 3 & 4 & 4 & 6  & 8  & 9 & 12 & 16 & 8 & 12 & 16 & 18 & 24 & 32 & 27 & 36 & 48 & 64\\
2 & &  &  &   & 2 & 6  & 12  & 18 & 36 & 72 & 16 & 42 & 80 & 108 & 204 & 384 & 270 & 504 & 936 & 1,728\\
3 & &  &  &   &  &   &   &  6 & 24 & 96 & 8 & 48 & 144 & 264 & 768 & 2,208 & 1,278 & 3,552 & 9,696 & 25,920\\
4 & &  &  &   &  &   &   &  & & 24 & 2 & 18 & 84 & 270 & 1,332 & 6,504 & 3,078 & 13,716 & 58,752 & 239,760\\
5 & &  &  &   &  &   &   &  & &   &   &   &   &  108 & 1,008 & 9,792 & 3,834 & 29,808 & 216,864 & 1,437,696\\
6 & &  &  &   &  &   &   &  & &   &   &   &   &  12 & 264 & 7,104 & 2,412 & 36,216 & 494,064 & 5,728,896\\
7 & &  &  &   &  &   &   &  & &   &   &   &   &    &   &  2,112 & 756 & 23,760 & 691,200 & 15,326,208\\
8 & &  &  &   &  &   &   &  & &   &   &   &   &    &   &  216 & 108 & 7,776 & 581,688 & 27,534,816\\
9 & &  &  &   &  &   &   &  & &   &   &   &   &    &   &    & 12 & 1,056 & 283,584 & 32,971,008\\
10 & &  &  &   &  &   &   &  & &   &   &   &   &    &   &    &   &   &  75,744 & 25,941,504\\
11 & &  &  &   &  &   &   &  & &   &   &   &   &    &   &    &   &   & 10,368 & 13,153,536\\
12 & &  &  &   &  &   &   &  & &   &   &   &   &    &   &    &   &   & 576 & 4,215,744\\
13 & &  &  &   &  &   &   &  & &   &   &   &   &    &   &    &   &   &   & 847,872\\
14 & &  &  &   &  &   &   &  & &   &   &   &   &    &   &    &   &   &   & 110,592\\
15 & &  &  &   &  &   &   &  & &   &   &   &   &    &   &    &   &   &   & 9,216\\
16 & &  &  &   &  &   &   &  & &   &   &   &   &    &   &    &   &   &   & 576\\ \hline
Total & 2 & 3 & 4
& 5 & 7 & 13 & 21 & 34 & 73 & 209 & 35 & 121 & 325 & 781 & 3,601 & 28,353 & 11,776 & 116,425 & 2,423,521 & 127,545,137\\ \hline
\end{tabular}
}
\end{center}
\vspace{-0.2cm}
\caption{Distribution of $\mathcal{R}_{r,s,n}$ according to the size, for $r\leq s\leq n\leq 4$.}
\label{table_hs1}
\end{table}

\section{Enumeration and classification of $\mathcal{SOR}_{r,n}$.}

The ideal $I_{r,r,n}$ of Theorem \ref{thm0} can be slightly modified to impose partial Latin rectangles to be self-orthogonal. In this regard, the proof of the next result is analogous to that of the mentioned theorem.

\begin{thm} \label{thm1} The set $\mathcal{SOR}_{r,n}$ can be identified with the set of zeros of the following zero-dimensional ideal of $\mathbb{F}_2[x_{111},\ldots,x_{rrn}]$.
$$I_{r,n}=I_{r,r,n}\cup \langle\,x_{ijp}\cdot x_{klp}\cdot x_{jiq}\cdot x_{lkq}\colon\, i,j,k,l\in [r],\, p,q\in [n],\,(i,j)\neq(k,l)\,\rangle.$$
Moreover, $|\mathcal{SOR}_{r,n}|=\mathrm{dim}_{\mathbb{F}_2} (\mathbb{F}_2[{\bf x}]/I_{r,n})$ and $|\mathcal{SOR}_{r,n:m}|=\mathrm{HF}_{\mathbb{F}_2[{\bf x}]/I_{r,n}}(m)$.

\hfill $\Box$
\end{thm}

\vspace{0.5cm}

\renewcommand{\baselinestretch}{1}
\renewcommand{\tabcolsep}{0.5pt}
\begin{table}[htb]
\begin{center}{\tiny
\begin{tabular}{rrrrrrrrrrrrrrrrrrrrrr} \hline
\ & \multicolumn{5}{l}{$|\mathcal{SOR}_{r,n:m}|$}\\ \cline{2-19}
\ & \multicolumn{5}{l}{$r.n$}\\ \cline{2-19}
$m$ & 2.1 & 2.2 & 2.3 & 2.4 & 2.5 & 2.6 & 2.7 & 2.8 & 2.9 & 3.1 & 3.2 & 3.3 & 3.4 & 3.5 & 3.6 & 3.7 & 3.8 & 3.9 \\ \hline
0 & 1 & 1 & 1 & 1 & 1 & 1 & 1  & 1 & 1 & 1 & 1 & 1 & 1 & 1 & 1 & 1 & 1 & 1\\
1 & 4 & 8 & 12 & 16 & 20 & 24 & 28  & 32 & 36 & 9 & 18 & 27 & 36 & 45 & 54 & 63 & 72 & 81\\
2 & & 12 & 36 & 72 & 120 & 180 & 252  & 336 & 432 & 12 & 96 & 252 & 480 & 780 & 1,152 & 1,596 & 2,112 & 2,700\\
3 & &  &  24 & 96 & 240 & 480  & 840  & 1,344 & 2,016 & 2 & 172 & 1,014 & 3,032 & 6,730 & 12,612 & 21,182 & 32,944 & 48,402 \\
4 & &  &  &  24 & 120 &  360 & 840  &  1,680 & 3,024 & & 108 & 1,836 & 9,720 & 31,320 & 77,220 & 161,028 & 299,376 & 511,920 \\
5 & & &  &  &  &  &  &   &   & & 12 & 1,476 & 15,912 & 80,040 & 270,900 & 720,972 & 1,633,296 & 3,296,592 \\
6 & & &  &  &  &  &  &   &   & & & 444 & 12,816 & 110,040 & 537,360 & 1,883,700 & 5,313,504 & 12,859,056 \\
7 & & &  &  &  &  &  &   &   & & & 36 & 4,608 & 76,680 & 573,120 & 2,743,020 & 9,870,336 & 29,142,288 \\
8 & & &  &  &  &  &  &   &   &  & & & 720 & 24,480 & 295,920 & 2,005,920 & 9,444,960 & 34,655,040 \\
9 & & &  &  &  &  &  &   &   &  & & & 48 & 3,120 & 58,320 & 566,160 & 3,551,520 & 16,456,608 \\
\hline
Total & 5 & 21 & 73 & 209 & 501 & 1,045 & 1,961 & 3,393 & 5,509 & 24 & 407 & 5,086 & 47,373 & 333,236 & 1,826,659 & 8,103,642 & 30,148,121 & 96,972,688 \\ \hline
\end{tabular}}
\end{center}
\caption{Distribution of $\mathcal{SOR}_{r,n}$ according to the size, for $r\leq 3$ and $n\leq 9$.}
\label{table0}
\end{table}

In order to compute the reduced Gr\"obner basis and the Hilbert series of the ideal of Theorem \ref{thm1}, we have implemented the procedure {\em SOR} in the library {\em pls.lib}. It has been used to obtain the distribution in Table \ref{table0} of the partial Latin rectangles of $\mathcal{SOR}_{r,n}$ according to their size, for $r\leq 3$ and $n\leq 9$. Although the CPU running time is similar to that of the analogous cases of Table \ref{table_hs1}, our computer system runs out of memory to determine the case $r=4$. It is due to the fact that the computation of reduced Gr\"obner bases and Hilbert series is extremely sensitive to the number of variables and polynomials. We propose two possible strategies to reduce this cost of computation.

\subsection{First strategy. Direct sum.}

Given $P\in\mathcal{R}_{r,r,n}$ and $Q\in\mathcal{R}_{r',r',n}$, let $P\oplus Q\in\mathcal{R}_{r+r',r+r',n}$ be the partial Latin rectangle having $P$ as upper left corner, $Q$ as lower right corner and empty the rest of its cells. Let us consider the set $S_{r,r',n}=\{(P,Q)
\in\mathcal{SOR}_{r,n}\times\mathcal{SOR}_{r',n}\colon\, P\oplus Q\in\mathcal{SOR}_{r+r',n}\}$ and the ideal $I_{r,r',n}^{P,Q}=I_{r,r',n} + \sum_{(i,j,k)\in [r]\times [r']\times [n] }\langle\, x_{ijk}\mid\, \exists\ l\in [r] \text{ such that } p_{il}=k \text{ or } \exists\ l\in [r'] \text{ such that } q_{lj}=k\rangle$. The affine variety $V(I_{r,r',n}^{P,Q})$ coincides with that of partial Latin rectangles of $\mathcal{R}_{r,r',n}$ which can be included in the upper right corner of $P\oplus Q$ to obtain an element of $\mathcal{R}_{r+r',r+r',n}$. Now, given $A=(a_{ij})\in V(I_{r,r',n}^{P,Q})$, let $I_{r',r,n}^{P,Q,A}$ be the following subideal of $I_{r',r,n}$.

\begin{center}
$I_{r',r,n} +
\langle\, x_{ijk}x_{i'j'k}\colon i,i'\in [r'], j,j'\in [r],(i,j)\neq (i',j'),k\in [n], a_{ji}=a_{j'i'}\in [n]\rangle\, + \langle\, x_{ijk}x_{i'j'k'}\colon i,i'\in [r'], j,j'\in [r], k,k'\in [n], a_{ji}=k'\neq k=a_{j'i'}\rangle\, +  \langle\, x_{ijk}\colon (i,j,k)\in [r']\times [r]\times [n],\, a_{ji}=k, \text{ or }\, \exists\ l,m\in[r] \text{ such that }a_{ji}=p_{lm},\, p_{ml}=k, \text{ or }\, \exists\ l,m\in[s] \text{ such that } a_{ji}=q_{lm},\, q_{ml}=k, \text{ or }\, \exists\ l\in [r] \text{ such that } p_{lj}=k, \text{ or }\, \exists\ l\in [r'] \text{ such that } q_{il}=k\rangle\,$
\end{center}
It is verified that
\begin{equation}\label{eqHS}
|\mathcal{SOR}_{r,n:m}|=\sum_{\small \substack{(P,Q)\in S_{r,r',n}\\ A\in V(I_{r,r',n}^{P,Q})}}
t^{|P|+|Q|+|A|}\cdot
\mathrm{HS}_{\mathbb{F}_2/I^{P,Q,R}_{r',r,n}}(t).
\end{equation}
This expression has been used to check the data of Table \ref{table0} and to expose in Table \ref{table0a} the cardinality of $\mathcal{SOR}_{4,n}$, for all $n\leq 9$.

\renewcommand{\baselinestretch}{1}
\renewcommand{\tabcolsep}{0.3pt}
\begin{table}[htb]
\begin{center}{\tiny
\begin{tabular}{rrrrrrrrrr}
\ & \multicolumn{5}{l}{$r.n$}\\ \cline{2-10}
$m$ & 4.1 & 4.2 & 4.3 & 4.4 & 4.5 & 4.6 & 4.7 & 4.8 & 4.9 \\ \hline
0 & 1 & 1 & 1 & 1 & 1 & 1 & 1 & 1 & 1 \\
1 & 16 & 32 & 48 & 64 & 80 & 96 & 112 & 128 & 144\\
2 & 60 & 360 & 900 & 1,680 & 2,700 & 3,960 & 5,460 & 7,200 & 9,180\\
3 & 56 & 1,792 & 8,568 & 23,744 & 50,680 & 92,736& 153,272 & 235,648 & 343,224\\
4 & 14 & 4,196 & 45,306 &  199,784 & 587,750 & 1,373,004 & 2,763,026 & 5,008,976 & 8,405,694\\
5 &  & 4,560 & 137,520 & 1,046,880 & 4,428,960 & 13,552,560 & 33,783,120 & 73,106,880 & 142,655,040\\
6 &  & 2,256 & 240,216 & 3,479,616 & 22,225,680 & 91,696,080 & 288,559,656 & 755,440,896 & 1,731,190,176\\
7 &  & 480 & 237,888 & 7,350,912 & 74,983,680 & 430,875,360 & 1,748,093,760 & 5,618,070,528 & 15,283,095,552\\
8 &  & 24 & 131,544 & 9,785,664 & 169,923,120 & 1,410,554,520 & 7,551,498,024 & 30,273,440,064 & 98,905,243,104\\
9 &  &  & 40,896 & 8,103,552 & 256,494,720 & 3,202,600,320& 23,211,048,000 & 118,117,015,296 & 469,324,461,312\\
10 &  &  & 7,056 & 4,147,584 & 254,539,680 & 4,988,125,440 & 50,312,927,280 & 331,193,485,056 & 1,622,312,241,984\\
11 &  &  & 576 & 1,332,864 & 163,762,560 & 5,241,536,640 & 75,710,577,600 & 657,677,857,536 & 4,029,212,001,024\\
12 &  &  &  48 & 283,200 & 67,632,480 & 3,633,984,960& 77,231,577,360 & 903,490,374,528 & 7,027,446,121,920\\
13 &  &  & & 43,008 & 17,850,240 &  1,613,064,960& 51,545,020,800 & 827,927,331,840 & 8,299,928,625,408\\
14 &  &  & &  5,760 & 2,975,040 & 437,253,120 & 21,258,498,240 & 476,757,469,440 & 6,249,614,071,680\\
15 &  &  & & 768 & 291,840 & 65,571,840 & 4,861,006,080 & 154,221,473,280 & 2,678,459,470,848\\
16 &  &  & & 48 & 14,160 & 4,127,760& 466,312,560 & 21,145,881,120 & 492,310,895,328 \\
\hline
Total & 147 & 13,701 & 850,567 & 35,805,129 & 1,035,763,371 & 21,134,413,357 & 314,221,824,351 & 3,527,256,198,417
 & 30,984,678,831,619\\ \hline
\end{tabular}}
\end{center}
\caption{Distribution of $\mathcal{SOR}_{4,n}$ according to the size, for $n\leq 9$.}
\label{table0a}
\end{table}

\subsection{Second strategy. Number of distinct symbols.}

Given $s\leq n$, let $\mathcal{SOR}_{r,n;s}$ be the subset of partial Latin rectangles of $\mathcal{SOR}_{r,n}$ which contain exactly $s$ distinct symbols in their cells. Observe that any partial Latin rectangle of $\mathcal{SOR}_{r,n;s}$ is equal, up to permutation of symbols, to a partial Latin rectangle of $\mathcal{SOR}_{r,s;s}$. If $s=0$, then the set $\mathcal{SOR}_{r,0;0}$ only contains the empty $r\times r$ Latin rectangle. We have that
\begin{equation}\label{eqSO}
\mathcal{SOR}_{r,n}=\bigcup_{s=0}^n\mathcal{SOR}_{r,n;s}
\end{equation}
Hence, if $\sigma_{r,s}$ denotes the cardinality of $\mathcal{SOR}_{r,s;s}$, then
\begin{equation}\label{eq1}
|\mathcal{SOR}_{r,n}|=\sum_{s=0}^n\binom{n}{s}\cdot \sigma_{r,s},
\end{equation}
The following result holds.

\begin{lmm}\label{lmm2} It is verified that:
\begin{enumerate}[i.]
\item $\sigma_{r,0}=1$.
\item $\sigma_{r,r^2-1}=\frac 12\cdot (r^2+1)! \cdot \frac {r^3-2r^2+r+2}r$.
\item $\sigma_{r,r^2}=r^2!$.
\item $\sigma_{r,s}=0$, for all $s>r^2$.
\end{enumerate}
\end{lmm}

\begin{pf} Claims (i), (iii) and (iv) follows straightforward from the definition of $\sigma_{r,s}$. To prove the second expression, observe that any $r\times r$ partial Latin rectangle which contains exactly $r^2-1$ distinct symbols has at most one empty cell. There are $r^2!$ possible partial Latin rectangles with exactly one empty cell. Otherwise, there exist two cells with the same symbol, where at most one of them is in the main diagonal. Specifically, there are $(r^2-1)!\cdot r\cdot (r-1)\cdot (r-2)$ possible partial Latin rectangles with one such a cell in the main diagonal and $\frac 12 \cdot (r^2-1)!\cdot (r^2-r)\cdot (r-1)\cdot (r-2)$ ones without any of them in the main diagonal. The result follows from the addition of all these possibilities.\hfill $\Box$
\end{pf}

\vspace{0.5cm}

The enumeration of $\mathcal{SOR}_{r,s;s}$ can be based on that of $\mathcal{SOR}_{r,s-1;s-1}$, for all $s\in [n]$. To see it, given a partial Latin rectangle $P=(p_{ij})\in\mathcal{SOR}_{r,s-1;s-1}$, let us define the subset $\mathcal{SOR}^P_{r,s;s}$ of partial Latin rectangles $Q=(q_{ij})\in \mathcal{SOR}_{r,s;s}$ such that $q_{ij}=p_{ij}$ if $p_{ij}\in [s-1]$ and $q_{ij}=\emptyset$, otherwise. It is then verified that
\begin{equation}\label{eq00}
\mathcal{SOR}_{r,s;s}=\bigcup_{P\in\mathcal{SOR}_{r,s-1;s-1}} \mathcal{SOR}^P_{r,s;s}.
\end{equation}
The partial Latin rectangles of $\mathcal{SOR}^P_{r,s;s}\cup \{P\}$ can be identified with the zeros of the ideal $I_{r,s;s}^P$ based on the $r^2-|P|$ variables which result after substituting the following variables in the ideal $I_{r,s}$ of Theorem \ref{thm1}.
\begin{enumerate}
\item[a)] Each variable $x_{ijk}$ such that $i,j\in [r]$ and $k\in [s-1]$ is substituted by $1$ if $p_{ij}=k$, or by $0$, otherwise.
\item[b)] Each variable $x_{ijs}$ such that $i,j\in [r]$ and $p_{ij}\in [s-1]$ is substituted by $0$.
\end{enumerate}
The reduced Gr\"obner basis of this new ideal can then be computed to determine explicitly the set $\mathcal{SOR}^P_{r,s;s}$. If the same reasoning is done for each partial Latin rectangle of $\mathcal{SOR}_{r,s-1;s-1}$ and for each $s\leq n$, then we can enumerate the set $\mathcal{SOR}_{r,n}$. The reduction on the number of variables of each ideal $I_{r,s;s}^P$, in comparison with the $r^2\cdot n$ variables of the ideal $I_{r,n}$, implies a significant improvement of the initial cost of computation which was necessary to enumerate the set $\mathcal{SOR}_{r,n}$ by applying Theorem \ref{thm1}. Even if such an improvement is obtained at the expense of time of computation (observe that it would be necessary to compute $\sum_{s=0}^{n-1}\sigma_{r,s}$ distinct reduced Gr\"obner bases), this time can be reduced if we make use of the distribution of $\mathcal{SOR}_{r,n}$ into main classes. In this regard, given $s\leq n$, let $\mathfrak{P}_{r,s;s}$ be the set of main classes of $\mathcal{SOR}_{r,s;s}$ and let $\mathfrak{P}_{s,P}$ denote the main class of $P\in\mathcal{SOR}_{r,s}$, where the necessity of the subscript $s$ for the number of symbols is due to the fact that $P$ is also a self-orthogonal partial Latin rectangle of $\mathcal{SOR}_{r,t}$, for all $t>s$. Since $S_r\times S_s\rtimes \overline{S}_3$ is a finite group which acts on $\mathcal{SOR}_{r,s}$, Burnside's lemma implies that
\begin{equation}\label{eq0}
\sigma_{r,s}=\sum_{P\in\mathfrak{P}_{r,s;s}}|\mathfrak{P}_{s,P}|= \sum_{P\in\mathfrak{P}_{r,s;s}} \frac {2\cdot r!\cdot s!} {|\mathfrak{I}_s(P,P)|+|\mathfrak{I}_s(P,P^t)|},
\end{equation}
where given two partial Latin rectangles $P,Q\in \mathcal{SOR}_{r,s}$, the set $\mathfrak{I}_s(P,Q)$ denotes the set of isotopisms which transform $P$ into $Q$. The next result shows how the polynomial method can be used to determine this set and its cardinality.

\begin{thm} \label{thm3} Given $P=(p_{ij}), Q=(q_{ij})\in\mathcal{SOR}_{r,s}$, the set $\mathfrak{I}_s(P,Q)$ can be identified with the set of zeros of the zero-dimensional ideal of $\mathbb{F}_2[x_{11},\ldots,x_{rr},$ $y_{11},\ldots,y_{ss}]$.
$$I_{\mathfrak{I}_s(P,Q)}=\langle\,1-\sum_{j\in [r]}x_{ij}\colon\,i\in [r]\,\rangle + \langle\,1-\sum_{j\in [n]}y_{ij}\colon\,i\in [s]\,\rangle + \langle\,1-\sum_{i\in [r]}x_{ij}\colon\,j\in [r]\,\rangle +$$ $$\langle\,1-\sum_{i\in [n]}y_{ij}\colon\,j\in [s]\,\rangle +\langle\,x_{ij}\cdot\left(1-x_{ij}\right)\colon\, i,j \in [r]\,\rangle +\langle\,y_{ij}\cdot\left(1-y_{ij}\right)\colon\, i,j \in [s]\,\rangle +$$ $$\langle\,x_{ik}\cdot x_{jl} \cdot (y_{p_{ij}q_{kl}}-1)\colon\, i,j,k,l\in [r], \text{ such that } p_{ij}, q_{kl}\in [s]\,\rangle +$$ $$\langle\,x_{ik}\cdot x_{jl} \colon\, i,j,k,l\in [r],\text{ such that } (p_{ij}=\emptyset \text{ and } q_{kl}\in [s]) \text{ or }$$ $$ \text{ or }(p_{ij}\in [s] \text{ and } q_{kl}=\emptyset)\,\rangle.$$
Moreover, $|\mathfrak{I}_s(P,Q)|= \mathrm{dim}_{\mathbb{F}_2}(\mathbb{F}_2[x_{11},\ldots,x_{rr},y_{11},\ldots,y_{ss}]/I_{\mathfrak{I}_s(P,Q)})$.
\end{thm}

\begin{pf}  Any isotopism $\Theta=(\alpha,\gamma)\in S_r\times S_s$ of $\mathcal{SOR}_{r,s}$ can be univocally identified with a zero $(x^{\Theta}_{11},\ldots,x^{\Theta}_{rr},y^{\Theta}_{11},\ldots, y^{\Theta}_{ss})$, where $x^{\Theta}_{ij}=1$, (respectively, $y^{\Theta}_{ij}=1$) if $\alpha(i)=j$ (respectively, $\gamma(i)=j$) and $0$, otherwise. The first six subideals of $I_{\mathfrak{I}_s(P,Q)}$ imply that $\Theta$ belongs to $S_r\times S_s$ and the following two subideals imply that $\Theta$ transforms $P$ into $Q$. The rest of the proof is analogous to that of Theorem \ref{thm0}. \hfill $\Box$
\end{pf}

\vspace{0.5cm}

It just remains to determine a representative partial Latin rectangle of each main class of $\mathcal{SOR}_{r,s;s}$. The following result shows how the enumeration of $\mathfrak{P}_{r,s-1;s-1}$ can determine that of $\mathfrak{P}_{r,s;s}$.

\begin{lmm}\label{lmm3} Let $\{P_1,\ldots,P_m\}$ be a set of representative partial Latin rectangles of $\mathfrak{P}_{r,s-1;s-1}$. Given $P\in \mathcal{SOR}_{r,s;s}$, there exist a natural $i\leq m$ and a partial Latin rectangle $Q\in \mathcal{SOR}_{r,s;s}^{P_i}$ such that $P$ is in the same main class than $Q$.
\end{lmm}

\begin{pf} Given $P\in \mathcal{SOR}_{r,s;s}$, let $P'\in \mathcal{SOR}_{r,s-1;s-1}$ be the self-orthogonal partial Latin rectangle which results after removing the symbol $s$ from the cells of $P$. There exist a natural $i\leq m$ and a paratopism $\Theta=(\alpha,\gamma,\pi)\in S_r\times S_{s-1}\rtimes \overline{S}_3$ such that $P_i=P'^{\Theta}$. It is enough to consider $Q=P^{\Theta'}\in\mathcal{SOR}_{r,s;s}^{P_i}$, where $\Theta'=(\alpha,\gamma',\pi)\in S_r\times S_s\rtimes \overline{S_3}$ is defined such that $\gamma'(i)=\gamma(i)$ if $i<s$ and $\gamma'(s)=s$.\hfill $\Box$
\end{pf}

\vspace{0.5cm}

Once we have determined a set $S\subseteq \mathcal{SOR}_{r,s-1;s-1}$ of representative partial Latin rectangles of  $\mathfrak{P}_{r,s-1;s-1}$, it is enough to enumerate the set $\bigcup_{P\in S}\mathcal{SOR}^P_{r,s;s}$ and to distribute its elements according to their main classes. Theorem \ref{thm3} can be used to determine this last distribution, because two partial Latin rectangles $P,Q\in\mathcal{SOR}_{r,s}$ are in the same main class if and only if $|\mathfrak{I}_s(P,Q)|+|\mathfrak{I}_s(P,Q^t)|>0$. Previously, in order to enumerate each set $\mathcal{SOR}^P_{r,s;s}$, we can make use of the ideal $I_{r,s;s}^P$ that was defined at the beginning of the current subsection. All these considerations have been implemented in the procedure {\em SOR} mentioned after Theorem \ref{thm1}. It has also been implemented a procedure called {\em ortisot} that computes the reduced Gr\"obner basis related to the ideal of Theorem \ref{thm3}. With the joint use of both procedures we have determined the main classes of $\mathcal{SOR}_{r,s;s}$, for all $r\leq 3$ and $s\leq r^2$. The number of these classes are shown in Table \ref{table1} and have been used to obtain from Expression (\ref{eq0}) the corresponding values of $\sigma_{r,s}$. The use of these values in Expression (\ref{eq1}) allow us to prove the next result.

\begin{thm} \label{thmGF} It is verified that:
\begin{enumerate}
\item $|\mathcal{SOR}_{1,n}|=n+1$.
\item $|\mathcal{SOR}_{2,n}|=n^4-2n^3+5n^2+1$.
\item $|\mathcal{SOR}_{3,n}|= n^9-15n^8+122n^7-604n^6+1973n^5-4201n^4+5640n^3-4240n^2+1347n+1$. \hfill $\Box$
\end{enumerate}
\end{thm}

\vspace{0.25cm}

\renewcommand{\baselinestretch}{1}
\renewcommand{\tabcolsep}{1pt}
\begin{table}[htb]
\begin{center}{\tiny
\begin{tabular}{l|rrr|rrr}
\ & \multicolumn{3}{l|}{$|\mathfrak{P}_{r,s;s}|$} & \multicolumn{3}{l}{$\sigma_{r,s}$}\\ \cline{2-7}
\ & \multicolumn{6}{l}{$r$} \\ \cline{2-7}
$s$ & 1 & 2 & 3 & 1 & 2 & 3\\ \hline
1 & 1 &   2 & 5 &  1! &   4 & 23\\
2 &  &   3 & 24 &  &   12 & 360\\
3 &  &  2 & 71 &  &  24 & 3,936\\
4 &  &  1 & 128 &  &  4! & 29,376\\
5 &  & & 122 &   & & 143,280\\
6 &  & & 67  &  & & 442,080\\
7 &  & & 22  &  & & 826,560\\
8 &  & & 4  &  & & 846,720\\
9 &  & & 1 &  & & 9!\\
\end{tabular}}
\end{center}
\caption{Main classes of $\mathcal{SOR}_{r,s;s}$.}
\label{table1}
\end{table}

\section{Conclusions.}

Keeping in mind the results of the current paper, we can conclude that the Combinatorial Nullstellensatz of Alon is a good method to deal with the enumeration and classification of self-orthogonal $r\times r$ partial Latin rectangles based on $n$ symbols. Based on the adjacency and independent sets of a Hamming graph, we have identified these combinatorial objects with the zeros of a zero-dimensional radical ideal. Besides, we have exposed two distinct strategies to reduce the cost of computation of the reduced Gr\"obner basis and Hilbert series of such an ideal. They can be used to enumerate explicitly the set of self-orthogonal partial Latin rectangles or to obtain some general formulas about its cardinality. All our results have been implemented into distinct procedures that have been applied in particular to determine the number of $r\times s$ partial Latin rectangles based on $n$ symbols, for $r,s,n\leq 4$ and that of self-orthogonal partial Latin rectangles of order $r\leq 4$.

\end{document}